\documentstyle[12pt,leqno,amssymb,amscd,graphics]{amsart}

\newtheorem{thm}{Theorem}[section]
\newtheorem{lemma}[thm]{Lemma}
\newtheorem{prop}[thm]{Proposition}

\theoremstyle{definition}
\newtheorem{dfn}[thm]{Definition}
\theoremstyle{remark}

\begin{document}

\newcommand{\ct}{\cite}
\newcommand{\pr}{\protect\ref}
\newcommand{\su}{\subseteq}
\newcommand{\pa}{{\partial}}
\newcommand{\im}{{Imm(F,\E)}}
\newcommand{\hf}{{1 \over 2}}
\newcommand{\Q}{{\mathbb Q}}
\newcommand{\R}{{\mathbb R}}
\newcommand{\Z}{{\mathbb Z}}
\newcommand{\G}{{\mathbb G}}
\newcommand{\F}{{\mathbb F}}
\newcommand{\cx}{{\mathbb C}}
\newcommand{\J}{{\mathbb J}}
\newcommand{\s}{{\mathbb S}}
\newcommand{\E}{{\mathbb E}}

\newcommand{\zo}{{\Z^o}}
\newcommand{\ze}{{\Z^e}}
\newcommand{\h}{\widehat}

\newcommand{\ca}{{a_1 \ch a_2}}
\newcommand{\cb}{{b_1 \ch b_2}}
\newcommand{\cc}{{c_1 \ch c_2}}

\newcommand{\da}{{a \ch k}}
\newcommand{\db}{{b \ch l}}

\newcommand{\br}{\overline}
\newcommand{\X}{{\mathbb X}}
\newcommand{\Y}{{\mathbb Y}}
\newcommand{\U}{{\mathbb U}}

\newcommand{\I}{{\mathrm{Id}}}
\newcommand{\4}{{\mathcal{H}}}
\newcommand{\C}{{\mathcal{C}}}

\newcommand{\1}{{(1)}}
\newcommand{\ch}{\choose}
\newcommand{\ce}{{\mathcal{C}^{\text{ev}}}}
\newcommand{\co}{{\mathcal{C}^{\text{od}}}}

\newcommand{\cce}{{\mathcal{C}^{\text{ev}}_0}}
\newcommand{\cco}{{\mathcal{C}^{\text{od}}_0}}

\newcommand{\xe}{{\X^{\text{ev}}}}
\newcommand{\xo}{{\X^{\text{od}}}}
\newcommand{\ye}{{\Y^{\text{ev}}}}
\newcommand{\yo}{{\Y^{\text{od}}}}
\newcommand{\ev}{{\text{ev}}}
\newcommand{\od}{{\text{od}}}

\newcommand{\we}{{\W^{\text{ev}}}}
\newcommand{\wo}{{\W^{\text{od}}}}

\newcommand{\hc}{{H_1(F,\C)}}
\newcommand{\ak}{{ \{ a_k \}  }}   
\newcommand{\bk}{{ \{ b_k \}  }}   
\newcommand{\tb}{{ \Leftrightarrow }} 
\newcommand{\bn}{{ \leftrightarrow }} 

\newcounter{numb}

\title{Order one invariants of planar curves}
\author{Tahl Nowik}
\address{Department of Mathematics, Bar-Ilan University, 
Ramat-Gan 52900, Israel}
\email{tahl@@math.biu.ac.il}
\date{December 25, 2007}
\urladdr{http://www.math.biu.ac.il/$\sim$tahl}

\begin{abstract}
We give a complete description of all order 1 invariants of planar curves. 
\end{abstract}

\maketitle

\section{Introduction}\label{intro}

A planar curve is an immersion of $S^1$ in $\R^2$.
The study of invariants of planar curves has been initiated by V. I. Arnold in \ct{a1},\ct{a2},
where he presented the three basic order 1 invariants $J^+,J^-,St$.
Various explicit formulas for these and other invariants appear in \ct{cd},\ct{p1},\ct{p3},\ct{s}.
The study of order 1 and higher order invariants has in general split into the study of $J$-invariants,
as in \ct{g1},\ct{g2},\ct{l},\ct{p2}, 
and $S$-invariants as in \ct{k},\ct{t},\ct{va}, 
where $J$-invariants are invariants which are unchanged when passing a triple
point, and $S$-invariants are invariants which are unchanged when passing a tangency.
But as we shall see, the space of invariants spanned by the order 1
$J$- and $S$-invariants is much smaller than the full space of order 1 invariants.

In this work we give a complete presentation of all order 1 invariants of planar curves.
Since the vector space in which the invariants take their values is arbitrary, the presentation
is in terms of a universal order 1 invariant.
The corresponding description for spherical curves, i.e. immersions of $S^1$ in $S^2$, appears in \ct{n6}.
Analogous study for immersions of surfaces in 3-space has been carried out in \ct{n1}--\ct{n5}.

The structure of the paper is as follows. In Section \pr{state} we 
present the basic definitions, and construct our invariant $F$ of planar curves.  
In Section \pr{abs} we construct an abstract invariant, as opposed to the explicit invariant $F$, 
and show that the abstract invariant is universal (Theorem \pr{absf}).
In Section \pr{prf} we prove our main result, namely, that a slight modification $\h{F}$ of
the explicit invariant $F$ is universal (Theorem \pr{main}). 
In Section \pr{ex} we prove the main algebraic fact used in the proof of Theorem \pr{main}.

\section{Definitions and statement of results}\label{state}
By a \emph{curve} we will always mean an immersion $c:S^1 \to \R^2$.   
Let $\C$ denote the space of all curves. 
A curve will be called \emph{stable} if its only self intersections are transverse double points.
The generic singularities a curve may have are either a tangency of first order between two strands, 
which will be called a $J$-type singularity, or three strands meeting at a point, 
each two of which are transverse, which will be called an $S$-type singularity.
Singularities of type $J$ and $S$ appear in Figures \pr{p2} and \pr{p3}. 
A generic singularity can be resolved in two ways, 
and there is a standard way for considering one resolution positive, and the other negative,
as defined in \ct{a1}.  
We denote by $\C_n \su \C$ ($n\geq 0$) the space of all curves which have precisely 
$n$ generic singularity points (the self intersection being elsewhere stable).
In particular, $\C_0$ is the space of all stable curves. 
An invariant of curves is a function $f:\C_0 \to W$,
which is constant on the connected components of $\C_0$, and where $W$ in this work 
will always be a vector space over $\Q$.

Given a curve $c\in \C_n$, with singularities located at $p_1,\dots,p_n \in \R^2$, 
and given a subset $A\su \{p_1,\dots,p_n\}$,
we define $c_A \in C_0$ to be the stable curve obtained from $c$ by resolving all singularities
of $c$ at points of $A$ into the 
negative side, and all singularities not in $A$ into the positive side.
Given an invariant $f:\C_0 \to W$ we define the ``$n$th derivative'' of $f$ to be the function
$f^{(n)}:\C_n \to W$ defined by
$$f^{(n)}(c)=\sum_{ A \su \{p_1,\dots,p_n\} } (-1)^{|A|} f(c_A)$$
where $|A|$ is the number of elements in $A$.
An invariant $f:\C_0\to W$ is called \emph{of order $n$} if 
$f^{(n+1)}(c)=0$ for all $c\in \C_{n+1}$.
The space of all $W$ valued invariants on $\C_0$ of order $n$ is denoted $V_n=V_n(W)$.
Clearly $V_n \su V_m$ for $n \leq m$.
In this work we give a full description of $V_1$. We will construct a ``universal'' 
order 1 invariant, by which we mean the following:

\begin{dfn}\label{uni}
An order 1 invariant $\h{f} :\C_0 \to \h{W}$ will be called \emph{universal}, if
for any $W$ and any order 1 invariant $f:\C_0 \to W$, there exists a unique linear map 
$\phi: \h{W} \to W$ such that $f = \phi \circ \h{f}$. In other words, for any $W$, the natural map
$Hom_\Q(\h{W},W) \to V_1(W)$ given by $\phi \mapsto \phi \circ \h{f}$, is an isomorphism.
\end{dfn}

\begin{dfn}\label{gen2}
A \emph{standard} annulus will be an annulus $A$ 
of the form $A = \R^2 - U$, where $U \su \R^2$ is a circular open disc 
(i.e. $U$ is a disc geometrically, not just topologically.)
An immersion $e:[0,1] \to A$ will be called a \emph{simple arc} if
\begin{enumerate}
\item $e(0),e(1) \in \pa A$.
\item The ends of $e$ are perpendicular to $\pa A$.
\item The self intersections of $e$ are transverse double points.
\end{enumerate}
\end{dfn}

Let $\zo \su \Z$ denote the set of odd integers. Let $D\Z$ denote the set of all columns ${a_1 \ch a_2}$
with $a_1 \in \Z$ and $a_2 \in \Z^o$.

\begin{dfn}\label{di}
Let $A$ be a standard annulus, and let $e$ be a simple arc in $A$.
\begin{enumerate}
\item For a double point $v$ of $e$
we define $i(v) \in \{1,-1\}$, where $i(v)=1$ if the orientation at $v$ given by the two tangents to $e$ at $v$,
in the order they are visited, coincides with the orientation of $A$ (restricted from $\R^2$). 
Otherwise $i(v)=-1$.
\item We define the \emph{top} index of $e$, $I_1(e)\in\Z$ by $I_1(e) = \sum_v i(v)$ 
where the sum is over all double points $v$ of $e$.
\item We define the \emph{bottom} index of $e$, $I_2(e)$, to be the index defined in \ct{k}, that is, 
$I_2(e) = (\phi - \omega)/\pi$ where $\phi$ is the total angle of rotation  of the radius vector 
(from the center of $U$ to $e$)
along $e$, and $\omega$ is the total angle of rotation of the tangent vector along $e$. 
Since by definition of simple arc the ends of $e$ are perpendicular to $\pa A$, 
we have that $I_2(e)$ is an odd integer, that is, $I_2(e) \in \zo$.
\item We define the \emph{double} index of $e$ to be $I(e) = {I_1(e) \ch I_2(e)} \in D\Z$.
\end{enumerate}
\end{dfn}

Define $\X$ to be the vector space over $\Q$ with basis all symbols $X^{a_1,b_1}_{a_2,b_2}$ where 
$a_1, b_1 \in \Z , a_2, b_2 \in \zo$.
For $c \in \C_0$, let $v$ be a double point of $c$, and let $u_1, u_2$ be the two tangents at  
$v$ ordered by the orientation of $\R^2$.
Let $U$ be a small circular disc neighborhood centered at $v$, 
and let $A=\R^2 - U$. Assume the strands of $c$ cross $\pa U$ perpendicularly, so
$c|_{c^{-1}(D)}$ defines two simple arcs $c_1,c_2$ in $A$. The ordering $c_1,c_2$ is chosen so that
the tangent $u_i$ leads to $c_i$, $i=1,2$. They will be called the \emph{exterior arcs} of $c$.
We denote $I(c_1)$ by $a_1(v) \ch a_2(v)$ and $I(c_2)$ by
$b_1(v) \ch b_2(v)$.
We define $F : \C_0 \to \X$ as follows:
$$F(c) = \sum_v X^{a_1(v),b_1(v)}_{a_2(v),b_2(v)}$$ 
where the sum is over all double points $v$ of $c$.

For $c \in \C_0$ let $\omega(c)$ denote the Whitney winding number of $c$,
and let $G : \C_0 \to \X$ be the order 0 invariant defined by $G(c)=X^{\omega(c),0}_{1,-1}$.
Let $\h{F} : \C_0 \to \X$ be defined by $\h{F}(c) = F(c) + G(c)$.
In this work we will prove that $\h{F}$ is a universal order one invariant.

Often one is interested in order 1 invariants only up to order 0 invariants, and in that case the addition of
$G$ is not needed. Our definition of universal order 1 invariant however requires that \emph{all}
order 1 invariants be obtained by composition with linear maps (in a unique way).
For this a correction such as $G$ is needed, as will be explained after the proof of Theorem \pr{main}.  
We will see (Lemma \pr{xw}) that if $F(c) \neq 0$ then
$\omega(c)$ can be computed from $F(c)$, so in terms of separating curves, the only difference between $F$ 
and $\h{F}$ is that $\h{F}$ separates the two different embedded curves, whereas $F$ vanishes on both.

\section{An abstract universal order one invariant}\label{abs}

In this section we will construct an ``abstract'' universal order 1 invariant, as opposed to the ``concrete''
invariant described in the previous section. 
This will be an intermediate step in proving that the concrete
invariant is universal.

\begin{figure}
\scalebox{0.8}{\includegraphics{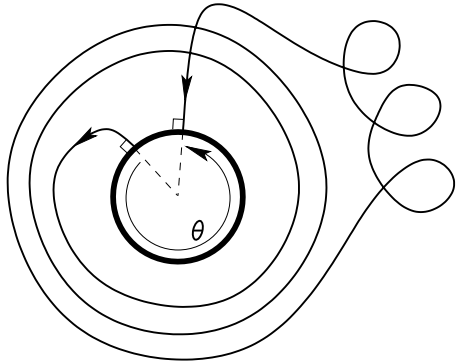}}
\caption{Representatives for the regular homotopy classes of simple arcs in $A$.}\label{p1}
\end{figure}

\begin{lemma}\label{omega}
Let $e$ be a simple arc in $A$ with initial point $i$ and final point $f$, 
with $i \neq f$.
Let $\theta$ be the angle from $i$ to $f$, 
in the positive direction, see Figure \pr{p1} . Let $\phi$ be the total angle of rotation  
of the radius vector (from the center of $U$ to $e$)
along $e$, and $\omega$ the total angle of rotation of the tangent vector along $e$.  
Then $I_1(e) = (2\phi - \omega - \theta + \pi)/2\pi$.
\end{lemma}

\begin{pf} Both sides of the equality are invariant under regular homotopy in $A$, 
keeping the initial and final points and tangents fixed (for the invariance of $I_1$ 
note that $i \neq f$). So, 
it is enough to check the equality on representatives of the various regular homotopy classes
of simple arcs, such as those of the form appearing in Figure \pr{p1}. 
They are arcs that begin by spiraling around $U$ a number of times 
and may then include a number of ``curls''. To verify the equality on such arcs, one can first
check it for embedded arcs, and then check that it is preserved 
when adding a round tour around $U$ and when adding a curl.
\end{pf}

\begin{lemma}\label{li}
If $e,e'$ are two simple arcs in $A$ with the same initial point $i$ and final point $f$, with $i \neq f$, 
then $I(e)=I(e')$ iff $e$ and $e'$ are regularly homotopic in $A$, keeping the initial and final
points and tangents fixed. 
\end{lemma}

\begin{pf}
Let $I_1=I_1(e),I_2=I_2(e)$. Let $\theta,\phi,\omega$ be as in Lemma \pr{omega}.
Then by Lemma \pr{omega}, $I_1 = (2\phi - \omega - \theta + \pi)/2\pi$, 
and by definition, $I_2 = (\phi - \omega)/\pi$.
These formulas are invertible, giving
$\phi = 2\pi I_1 - \pi I_2 + \theta - \pi$
and $\omega= 2\pi I_1 - 2\pi I_2 + \theta - \pi$.
Since the pair $\phi,\omega$ characterizes the regular homotopy class of $e$, so does the pair
$I_1,I_2$.
\end{pf}

Two curves $c,c' \in \C_1$ will be called \emph{equivalent} 
if there is an ambient isotopy of $\R^2$ bringing a neighborhood $U$ of 
the singular point of $c$ onto a neighborhood $U'$ of the singular point of $c'$, such that the configuration
near the singular point precisely matches and such that the exterior arcs in $A = \R^2-U$ are 
regularly homotopic in $A$. 
We always assume, as before, 
that the strands of $c$ cross $\pa U$ perpendicularly, and so the exterior arcs are simple arcs.
Note that the crossing points with $\pa U$ are all distinct (even for $J$ type singularities).  
By Lemma \pr{li} a pair of exterior arcs with same initial and final points, are regularly homotopic
in $A$ iff their double indices are the same. So we have:

\begin{prop}\label{eq1}
Two curves $c,c' \in \C_1$ are equivalent
iff  they have the same singularity configuration and the same corresponding double indices of exterior arcs.
\end{prop}

The following is clear from the definition of order 1 invariant:

\begin{lemma}\label{eqv}
Let $f : \C_0 \to W$ be an invariant, then $f$ is of order 1 iff for any two equivalent $c,c' \in \C_1$,
$f^{(1)}(c) = f^{(1)}(c')$.
\end{lemma}

We will attach a symbol to each equivalence class of curves in $\C_1$ as follows. For $J$ type singularities there are 
three distinct configurations which we name $J^+,J^A,J^B$, see Figure \pr{p2}.
The $J^A,J^B$ singularities are symmetric with respect to a $\pi$ rotation, which interchanges the two exterior arcs,
and so, by Proposition \pr{eq1},
the equivalence class of a $J^A$ or $J^B$ singularity is characterized by a symbol $J^A_{\ca,\cb}$ and $J^B_{\ca,\cb}$
($\ca,\cb \in D\Z$), where $\ca,\cb$ is an \emph{unordered} pair, registering the double indices of the two exterior arcs.
The $J^+$ configuration, on the other hand, is not symmetric, and so characterized by a symbol $J^+_{\ca,\cb}$ with 
$\ca,\cb$ an \emph{ordered} pair, with say, $\ca$ corresponding to the lower strand in Figure \pr{p2}.
(We will shortly see that this ordering may however be disregarded.)

As to $S$ type singularities, there are four types, as seen in Figure \pr{p3}, 
depending on the relative orientations of the three strands participating in
the triple point, and the way they are connected 
to each other. The distinction between them will be incorporated into the
way we register the double indices of the exterior arcs. 
We will have a cyclicly ordered triple
$\ca,\cb,\cc \in D\Z$ registering the double indices of the three exterior arcs, in the cyclic
order they appear along $S^1$, and each may appear with or without a hat, according to the following
rule. For given exterior arc $e$ with double index $\ca$, 
let $u_1$ be the initial tangent of $e$ and $u_2$ the initial tangent
of the following segment. Then if $u_2$ is pointing to the right of $u_1$ then $\ca$ will appear with a hat,
and if $u_2$ is pointing to the left of $u_1$ then $\ca$ will appear unhatted. Since the ordering is cyclic, this gives
four types of $S$ symbols, 
$S_{\ca,\cb,\cc} , S_{\h{\ca},\cb,\cc} , S_{\h{\ca},\h{\cb},\cc}, S_{\h{\ca},\h{\cb},\h{\cc}}$, which
correspond to the four types of $S$ singularities.

\begin{figure}
\scalebox{0.8}{\includegraphics{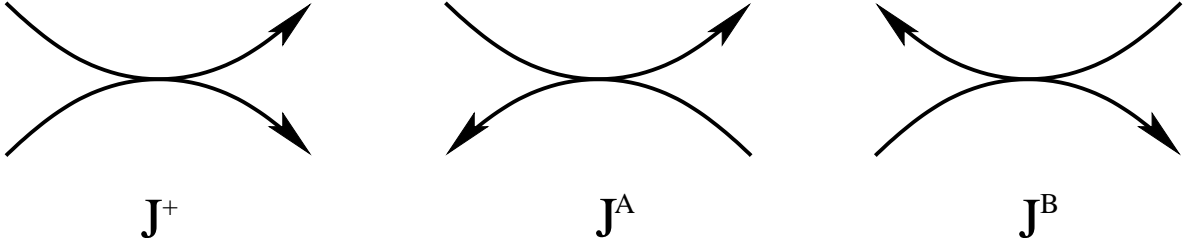}}
\caption{Singularities of type $J^+$, $J^A$, $J^B$.}\label{p2}
\end{figure}

\begin{figure}
\scalebox{0.8}{\includegraphics{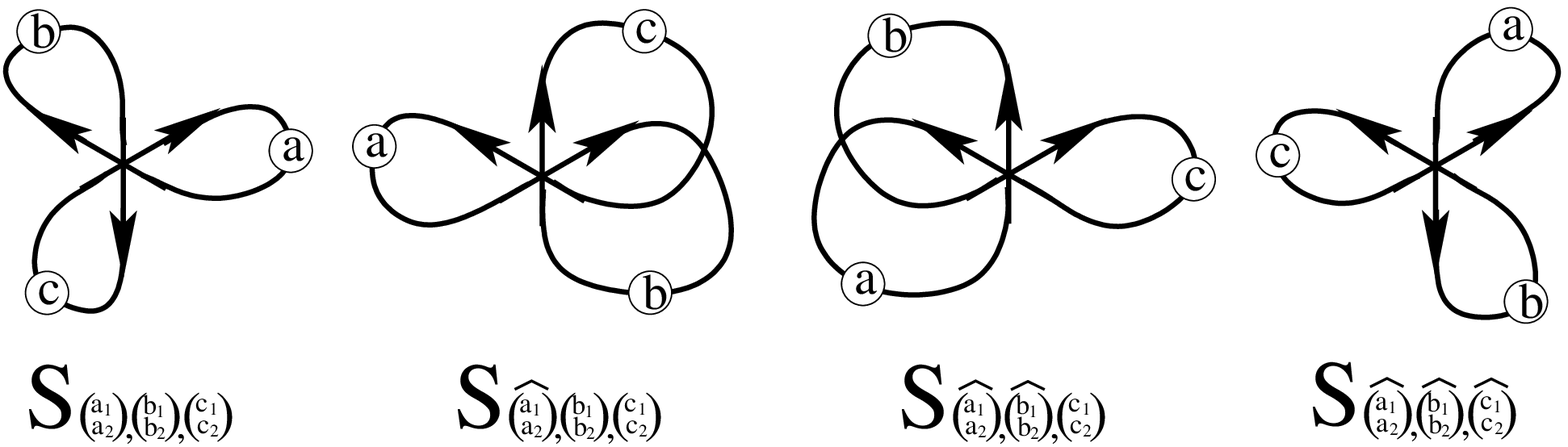}}
\caption{Singularities of type $S_{\ca,\cb,\cc}$, $S_{\h{\ca},\cb,\cc}$, $S_{\h{\ca},\h{\cb},\cc}$, 
$S_{\h{\ca},\h{\cb},\h{\cc}}$.}\label{p3}
\end{figure}

Let $\F$ denote the vector space over $\Q$ with basis all the above symbols, $J^+_{\ca,\cb}$, $J^A_{\ca,\cb}$, 
$J^B_{\ca,\cb}$, 
$S_{\ca,\cb,\cc}$, $S_{\h{\ca},\cb,\cc}$, $S_{\h{\ca},\h{\cb},\cc}$, $S_{\h{\ca},\h{\cb},\h{\cc}}$.
Let $\gamma$ be a generic path in $\C$, that is, 
a path (regular homotopy) $[0,1] \to \C$, whose image lies in $\C_0 \cup \C_1$ and 
which is transverse with respect to $\C_1$.
We denote by $v(\gamma) \in \F$ the sum of symbols of the singularities $\gamma$ passes, 
each added with $+$ or $-$
sign according to whether we pass it from its negative side to its positive side, or from its positive side to its 
negative side, respectively. Let $N \su \F$ be the subspace generated by all elements $v(\gamma)$ obtained from 
all possible generic \emph{loops} $\gamma$ in $\C$ (i.e. \emph{closed} paths), and let $\G = \F / N$. 

For an order 1 invariant $f:\C_0 \to W$, since by Lemma \pr{eqv} $f^{(1)}$ coincides on equivalent curves in $\C_1$, 
it induces a well defined linear map $f^{(1)}:\F \to W$. 

The following is clear:
\begin{lemma}\label{vgm}
If $\gamma$ is a generic path in $\C$, from $c_1$ to $c_2$, then $f^\1 ( v(\gamma) ) = f(c_2) - f (c_1)$. 
\end{lemma}

From Lemma \pr{vgm} it follows that $f^\1$ vanishes on the generators of $N$,
so it also induces a well defined linear map $f^{(1)} : \G \to W$.

Let $\Gamma_m \in \C_0$ be a curve with Whitney number $m$, chosen once and for all as
base curve for its regular homotopy class. 
Let $\h{\G} = \G \oplus \bigoplus_{m \in \Z} \Q u_m$,
where $\{ u_m \}_{m \in \Z}$ is a set of new independent vectors. 
We define an order 1 invariant $\h{f} : \C_0 \to \h{\G}$ as follows:
For any $c \in \C_0$ there is a generic path $\gamma$ 
from $\Gamma_m$ to $c$, where $m$ is the Whitney number of $c$. 
Let $\h{f}(c)=u_m + v(\gamma) \in \h{\G}$. 
By definition of $N$, 
$v(\gamma) \in \G$ is indeed independent of the choice of path $\gamma$.
From Lemma \pr{eqv} it is clear that $\h{f}$
is an order 1 invariant, and we will now see that it is \emph{universal}: 

\begin{thm}\label{absf}
$\h{f} : \C_0 \to \h{\G}$ is a universal order 1 invariant. 
\end{thm}

\begin{pf}
For an order 1 invariant $f: \C_0 \to W$, define the linear map $\Phi_f:\h{\G}\to W$
by $\Phi_f|_{\G} = f^{(1)}$ and $\Phi_f(u_m)=f(\Gamma_m)$, $m \in \Z$.
We claim $\Phi_f \circ \h{f} = f$ and that $\Phi_f$ is the unique linear map satisfying this property.
Indeed, let $c \in \C_0$ and let $\gamma$ be a generic path from $\Gamma_m$ to $c$ where $m$ is the 
Whitney number of $c$.
Then by Lemma \pr{vgm} 
$f(c) = f(\Gamma_m) + f^\1( v(\gamma) ) = \Phi_f(u_m) + \Phi_f ( v(\gamma)) = \Phi_f ( \h{f} (c))$. 

For uniqueness, it is enough to show that $\h{f} ( \C_0 )$ spans $\h{\G}$.
We have $u_m = \h{f}(\Gamma_m)$ for  $m \in \Z$ so it remains to show that
$\G \su span \h{f} ( \C_0 )$.
Indeed for any generating symbol $T$ of $\G$, $T$ is the 
the difference $\h{f}(c)-\h{f}(c')$ for two
curves $c,c' \in \C_0$, namely, the two resolutions of a curve in $\C_1$ 
whose symbol is $T$.
\end{pf}

The following is also clear:

\begin{lemma}\label{iso}
Let $f: \C_0 \to W$ be an order 1 invariant.
Then $f$ is a universal order 1 invariant iff
$\Phi_f: \h{\G} \to W$ (appearing in the proof of Theorem \pr{absf}) is an isomorphism. 
\end{lemma}

We present two specific subfamilies of the set of generators of $N$.
They appear in Figures \pr{p4},\pr{p5}. 
The figures are schematic, only indicating the way the strands 
involved in the local configuration are 
connected to each other. The two exterior arcs for the intersection point appearing in the left most curves,
labeled $a$ and $b$, are assumed to have double indices $\ca$ and $\cb$ respectively.
To obtain the double indices for the tangency points, note that exterior arcs that may be obtained from each
other by sliding their ends along $\pa A$ have the same bottom index.
The element obtained from Figure \pr{p4} is $J^+_{\ca,\cb} - J^+_{\cb,\ca}$, 
that is, in $\G$ we have the relation $J^+_{\ca,\cb} - J^+_{\cb,\ca} = 0$,
and so from now on we will simply regard the double indices of $J^+_{\ca,\cb}$ as being \emph{un-ordered} 
(as is true by definition for the $J^A$ and $J^B$ symbols). 
The relation obtained from Figure \pr{p5} is $J^B_{{a_1+1 \ch a_2},\cb} - J^A_{\ca,{b_1-1 \ch b_2}} = 0$, 
or after index shift:
$J^B_{\ca,\cb} = J^A_{{a_1-1 \ch a_2},{b_1-1 \ch b_2}}$.

\begin{figure}
\scalebox{0.8}{\includegraphics{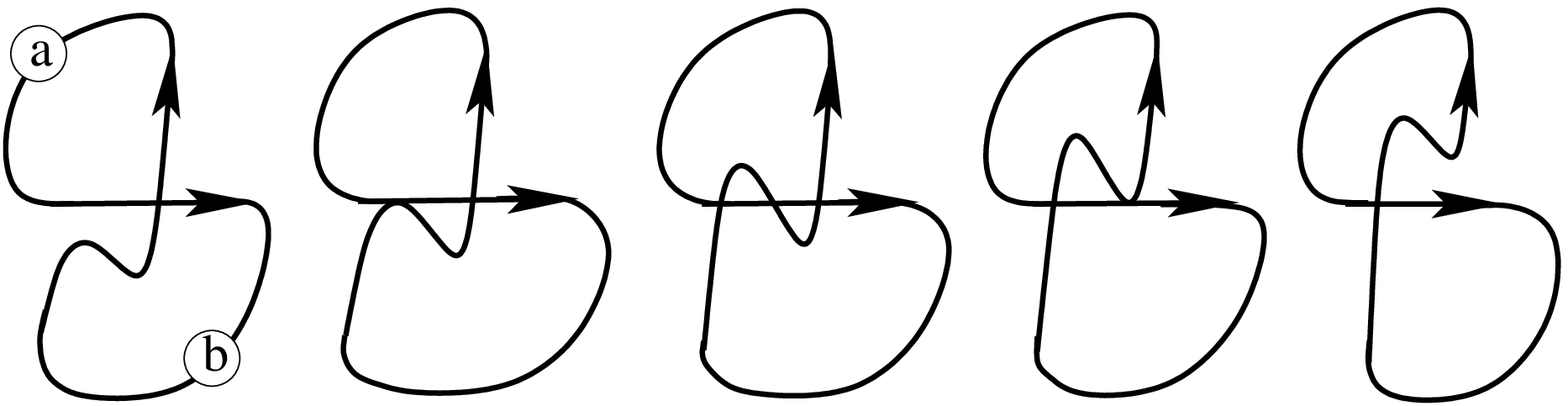}}
\caption{$J^+_{\ca,\cb} - J^+_{\cb,\ca} = 0$}\label{p4}
\end{figure}

\begin{figure}
\scalebox{0.8}{\includegraphics{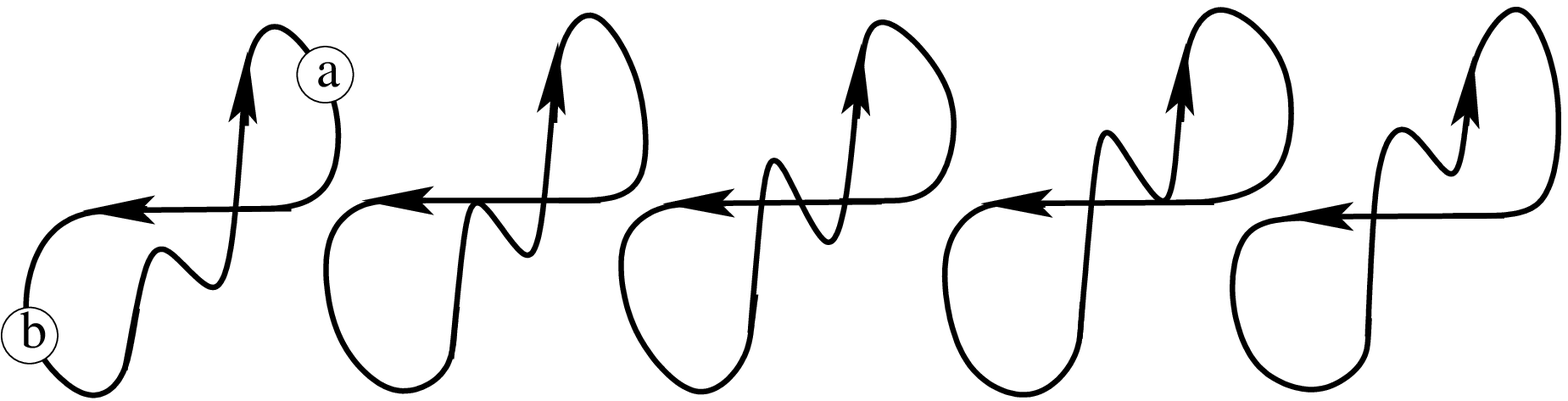}}
\caption{$J^B_{{a_1+1 \ch a_2},\cb} - J^A_{\ca,{b_1-1 \ch b_2}} = 0$}\label{p5}
\end{figure}

\section{Proof of main result}\label{prf}

Returning to our invariant $F:\C_0 \to \X$, we look at $F^{(1)}$.
It is easy to see that indeed the value of $F^\1$ depends only on the symbol of a
given curve $c \in \C_1$, which by Lemma \pr{eqv} proves  
that indeed $F$ is an order 1 invariant. 
We will need the value of $F^\1$ only on the symbols of type $J^+$ and $J^A$, which are demonstrated in 
Figures \pr{p6},\pr{p7}. Again, the figures are schematic, only indicating the way the strands 
involved in the local configuration
are connected to each other. The two exterior arcs for the tangency point appearing in the middle curves,
labeled $a$ and $b$, are assumed to have double indices $\ca$ and $\cb$ respectively.
In Figure \pr{p6} the two new intersection points $u$ and $v$ contribute additional terms 
$X^{a_1,b_1}_{a_2,b_2}$ and $X^{b_1,a_1}_{b_2,a_2} $ respectively, and in Figure \pr{p7} $u$ and $v$ contribute
$X^{a_1,b_1+1}_{a_2,b_2}$ and $X^{b_1,a_1+1}_{b_2,a_2}$ respectively. 
So we have:  
\begin{itemize}
\item $F^\1(J^+_{{a_1 \ch a_2},{b_1 \ch b_2}}) = 
X^{a_1,b_1}_{a_2,b_2} + X^{b_1,a_1}_{b_2,a_2} $
\item $F^\1(J^A_{{a_1 \ch a_2},{b_1 \ch b_2}}) = 
X^{a_1,b_1+1}_{a_2,b_2} + X^{b_1,a_1+1}_{b_2,a_2} $
\end{itemize}

\begin{figure}
\scalebox{0.8}{\includegraphics{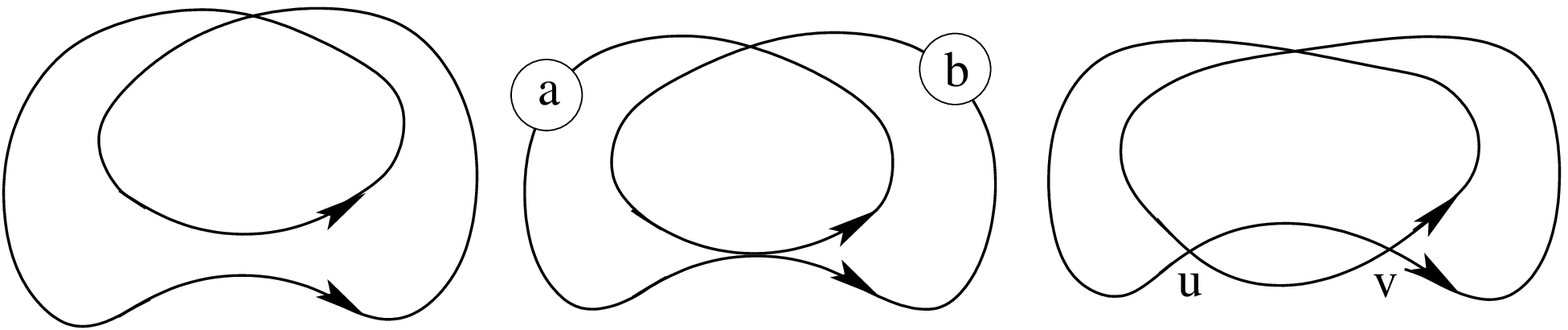}}
\caption{$F^\1(J^+_{{a_1 \ch a_2},{b_1 \ch b_2}}) = 
X^{a_1,b_1}_{a_2,b_2} + X^{b_1,a_1}_{b_2,a_2} $}\label{p6}
\end{figure}

\begin{figure}
\scalebox{0.8}{\includegraphics{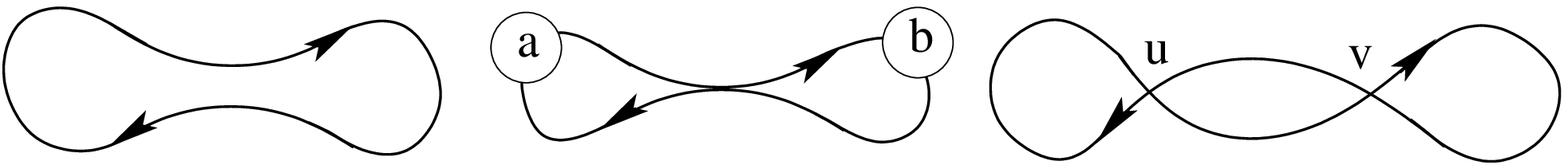}}
\caption{$F^\1(J^A_{{a_1 \ch a_2},{b_1 \ch b_2}}) = 
X^{a_1,b_1+1}_{a_2,b_2} + X^{b_1,a_1+1}_{b_2,a_2} $}\label{p7}
\end{figure}

We define $\Y$ to be the vector space over $\Q$ generated by the symbols $Y^n_{k,l}$, $n \in \Z$, 
$k,l \in \zo$ and with the relation $$Y^n_{k,l} = -Y^n_{l,k}.$$
We define the map $\Psi: \X \to \Y$ by $$\Psi(X^{a,b}_{k,l}) = Y^{a+b}_{k,l}.$$

For each $m \in \Z$ define the linear functional $g_m :\Y \to \Q$ as follows:
$$g_m(Y^n_{k,l}) =
\begin{cases}
k-l         & \ \ \ \  \text{if} \ \ \ n-k-l=m \ \ \ \text{and} \ \ \ k+l =0       \\
(l-k)/2     & \ \ \ \  \text{if} \ \ \ n-k-l=m \ \ \ \text{and}  \ \ \ k+l = \pm 2       \\
0           & \ \ \ \  \text{otherwise}   
\end{cases}$$
Let $\br{\Y} \su \Y$ denote the subspace defined by the family of all equation $g_m = 0$, $m \in \Z$,
that is, the intersection of kernels of all $g_m$.
Let $\br{\X} = \Psi^{-1}(\br{\Y})$.

We now look at the order 1 invariant $K = \Psi \circ F : \C_0 \to \Y$. 
The invariant $K$ is closely related to the invariant defined by M. \`{E}. Kazarian in \ct{k}.
Let $\J$ denote the subspace of $\G$ spanned by all symbols of type $J$, then by the relation
$J^B_{\ca,\cb} = J^A_{{a_1-1 \ch a_2},{b_1-1 \ch b_2}}$
it is also spanned by all symbols of type $J^+,J^A$.
It is clear from our explicit formulas for $F^\1$ on the symbols of type $J^+,J^A$, that
$K^\1 = (\Psi \circ F)^\1 = \Psi \circ F^\1$ vanishes on $\J$, and so induces a map 
$K^\1 : \G / \J \to \Y$. As we shall now explain, the result 
in \ct{k} may be stated, in the language of the present work, as follows:

\begin{thm}[M. \`{E}. Kazarian]\label{kaz}
$K^\1 : \G / \J \to \Y$ is injective, with image $\br{\Y}$.
\end{thm}

To make precise the relation between our $K$ and the invariant of \ct{k} we first prove:

\begin{lemma}\label{xw}
If $c \in \C_0$ is a curve with Whitney winding number $m$, then for each intersection point $v$ of $c$, 
the two double indices $\ca,\cb$ attached to $v$ satisfy $a_1+b_1-a_2-b_2=m$.
\end{lemma}

\begin{pf}
For the two exterior arcs $c_1,c_2$ of $v$, let $\theta_1,\phi_1,\omega_1$ and $\theta_2,\phi_2,\omega_2$
be the corresponding angles as in Lemma \pr{omega},
then $\theta_1+\theta_2=2\pi$ and $(\omega_1+\omega_2)/2\pi=m$. 
So $I_1(c_1) + I_1(c_2) - I_2(c_1) - I_2(c_2) = 
(2\phi_1 - \omega_1 - \theta_1 + \pi)/2\pi + (2\phi_2 - \omega_2 - \theta_2 + \pi)/2\pi
-(\phi_1 - \omega_1)/\pi -(\phi_2 - \omega_2)/\pi = m$.
\end{pf}

Accordingly, we define $\X^m \su \X$ to be the subspace spanned by all $X^{a,b}_{k,l}$ with $a+b-k-l=m$ and 
$\Y^m \su \Y$ to be the subspace spanned by all $Y^n_{k,l}$ with $n-k-l=m$.  
We have $\X=\oplus_m \X^m$, $\Y = \oplus_m \Y^m$, $F$ maps all curves with
Whitney number $m$ into $\X^m$ and $\Psi$ maps $\X^m$ into $\Y^m$. 
In a similar way each generating symbol
of $\G$ is related to a specific regular homotopy class and accordingly 
$\G$ and $\J$ split into $\oplus_m \G^m$ and $\oplus_m \J^m$.
Finally, define $\br{\Y}^m \su \Y^m$ to be the subspace defined by the equation $g_m=0$,
then $\br{\Y} = \oplus_m \br{\Y}^m$.
Each $\Y^m$ is naturally isomorphic to  $\Lambda^2 A$ of \ct{k} Remark 1, 
where $Y^{m+k+l}_{k,l} \in \Y^m$ is identified with $e_k \wedge e_l \in \Lambda^2 A$. 
This induces a projection $p: \Y \to \Lambda^2 A$, mapping any $Y^n_{k,l}$ to $e_k \wedge e_l$. 
The invariant $S$ of \ct{k} Remark 1 is precisely $p \circ K$. 
Now \ct{k} Theorem 2 relates to each regular homotopy class separately,
the space $I_1/I_0$ appearing there is the space of $S$-invariant on a given regular homotopy class, up to constants.
An $S$-invariant is by definition an invariant $f$ such that $f^\1$ vanishes on $\J$, and so $I_1/I_0$ is
precisely the dual to our $\G^m / \J^m$. So, the statement of \ct{k} Theorem 2 is that
the dual to our $K^\1 : \G^m / \J^m \to \br{\Y}^m$ is an isomorphism. It follows that 
this map itself is an isomorphism, and so $K^\1 : \G / \J \to \br{\Y}$ is an isomorphism.

We now look at the following commutative diagram:

$$\begin{CD}
0 @>>> \J @>inc>> \G @>quot>> \G / \J @>>> 0 \\
& &      @V Id VV     @VV F^\1 V       @VV K^\1 V \\    
0 @>>> \J @>>F^\1|_{\J}> \X @>>\Psi> \Y @>>> 0 
\end{CD}$$

The algebraic portion of our work is the following theorem, whose proof we defer to Section \pr{ex}:

\begin{thm}\label{exact}
The sequence \ \ $0 \to \J @>F^\1|_{\J}>> \X @>\Psi>> \Y \to 0$ \ \ is exact.
\end{thm}

It follows from Theorem \pr{kaz}, Theorem \pr{exact}, and the commutative diagram, that:

\begin{thm}\label{f1}
$F^\1 : \G \to \X$ is injective, with image $\br{\X}$.
\end{thm}

We may now prove our main result: 

\begin{thm}\label{main}
For $c \in \C_0$ let $\omega(c)$ denote the Whitney winding number of $c$,
and let $G : \C_0 \to \X$ be the order 0 invariant defined by $G(c)=X^{\omega(c),0}_{1,-1}$.
Let $\h{F} : \C_0 \to \X$ be defined by $\h{F}(c) = F(c) + G(c)$, then $\h{F}$ 
is a universal order one invariant.
\end{thm}

\begin{pf}
By Lemma \pr{iso} it is enough to show that $\Phi_{\h{F}} : \h{\G} \to \X$ is an isomorphism. 
Indeed, by Theorem \pr{f1}, 
$\Phi_{\h{F}}|_{\G} = \h{F}^\1 = F^\1$ is injective, and onto $\br{\X}$.
Now, the values $X^{\omega(c),0}_{1,-1}$ for $G$ were precisely chosen so that 
$\{ \h{F}(\Gamma_m) \}_{m \in \Z} = \{ \Phi_{\h{F}}(u_m) \}_{m \in \Z}$  will be a basis  for $\X / \br{\X}$.
This is verified by evaluating the functionals $g_m \circ \Psi$ 
(which are the defining equations of $\br{\X}$ in $\X$)
on $\h{F}(\Gamma_{m'})$. Indeed, from \ct{k} 
and our definition of $g_m$ we know
$g_m(\Psi( F(\Gamma_{m'}))) = g_m(K(\Gamma_{m'} ))= 2\delta_{m,m'}\delta_{m',0}$
and so $g_m \circ \Psi( \h{F}(\Gamma_{m'}) ) =  2\delta_{m,m'}\delta_{m',0} + 2\delta_{m,m'} 
= 2(1+\delta_{m',0})\delta_{m,m'}$.
It follows that $\Phi_{\h{F}} : \h{\G} \to \X$ is an isomorphism. 
\end{pf}

We now see why we needed to add to $F$ a correcting order 0 invariant such as $G$. 
From the proof of Theorem \pr{main} we see that $\Phi_F : \h{\G} \to \X$ 
is non-injective and non-surjective.
$\Phi_F$ being non-injective means that not all order 1 invariants may be realized as $\phi \circ F$ 
($\phi : \X \to W$ a linear map).
Indeed $g_m \circ \Psi( F(\Gamma_{m'})) = 2\delta_{m,m'}\delta_{m',0}$ implies that any order 0 invariant of
the form $\phi \circ F$ vanishes on all curves with Whitney number $m \neq 0$. 
$\Phi_F$ being non-surjective means that the choice of $\phi$ is non-unique, that is, there are $\phi \neq \psi$ with
$\phi \circ F = \psi \circ F$. Adding $G$ to $F$ solved both injectivity and surjectivity problems.
(If non-surjectivity were the only problem, then it would be fixed simply by diminishing the target space $\X$.) 

We conclude this section by explaining the remark made in the introduction, that the space of
invariants spanned by all order 1 $J$- and $S$-invariants is much smaller than the full space of order 1 invariants.
Indeed, let $\s \su \G$ be the subspace spanned  by all symbols of type $S$. By definition, an order 1 invariant $f$
is a $J$- or $S$-invariant, iff $f^\1$ vanishes on $\s$ or $\J$ respectively. 
Therefore an invariant $f$ is the sum of $J$- and $S$-invariants iff $f^\1$ vanishes on $\s \cap \J$.
By observing relations in $N$ coming from loops as in \ct{n6} Figures 6,7,
one can show that $\s \cap \J$ is the subspace of $\G$ 
spanned by all symbols of the form 
$J^+_{\ca,\cb} - J^+_{{a_1+2 \ch a_2+2},\cb}$ and $J^A_{\ca,\cb} -J^A_{{a_1+2 \ch a_2+2},\cb}$, 
$a_1, b_1 \in \Z , a_2, b_2 \in \zo$.
Given that the symbols of type $J^+,J^A$ are independent in $\G$ (a fact that will be proved in Section \pr{ex}),
we see $\s \cap \J$ is a fairly large subspace of $\G$.

\section{The exact sequence}\label{ex}

In this section we prove Theorem \pr{exact}, stating that the sequence 
$0 \to \J @>F^\1|_{\J}>> \X @>\Psi>> \Y \to 0$ is exact.
For $n \in \Z$ and $k \leq l \in \zo$ 
let $A^n_{k,l}$ be the following set of symbols, viewed as elements in $\G$. 
$$A^n_{k,l} = \{ J^+_{ {a \ch k} , {b \ch l} } \ : \  a+b = n \} \cup 
\{ J^A_{{a \ch k},{b \ch l}} \ : \ a+b=n-1 \}.$$   
Let $\J^n_{k,l} = span A^n_{k,l} \su \J$. We point out that since the pair ${a \ch k} , {b \ch l}$ is unordered,
the meaning of the pair $a,b$ in the conditions $a+b=n$ and $a+b=n-1$ depends on whether $k=l$ or $k<l$.
If $k=l$ then $a,b$ runs over all unordered pairs $\{a,b\}$, whereas if $k<l$ then the pairing of $a$ with $k$
and $b$ with $l$ serves to distinguish them and so in this case $a,b$ runs over all \emph{ordered} pairs $( a,b )$.

For $n \in \Z$ and $k \leq l \in \zo$ define $\X^n_{k,l}$
to be the span of all symbols of the form  $X^{a,b}_{k,l}$ and $X^{a,b}_{l,k}$ with $a+b =n$.
Then $\X$ is the direct sum of all $\X^n_{k,l}$.
Here there is also a difference between the cases $k=l$ and $k<l$. When $k=l$ a basis for $\X^n_{k,k}$ is 
$\{ X^{i,n-i}_{k,k} \}_{i \in \Z}$, whereas when $k < l$ a basis for $\X^n_{k,l}$ 
is $\{ X^{i,n-i}_{k,l} \}_{i \in \Z} \cup \{ X^{i,n-i}_{l,k} \}_{i \in \Z}$.
Let $B^n_{k,l} = F^\1(A^n_{k,l}) \su \X$. That is
$$B^n_{k,l} = \{ X^{a,b}_{k,l} + X^{b,a}_{l,k} \ : \ a+b=n \} 
\cup \{ X^{a,b+1}_{k,l} + X^{b,a+1}_{l,k} \ : \ a+b=n-1 \}.$$
We see  $B^n_{k,l} \su \X^n_{k,l}$, and so $F^\1(\J^n_{k,l}) \su \X^n_{k,l}$.
$B^n_{k,l}$ being simply the image of $A^n_{k,l}$, the meaning of the pair $a,b$ here is 
the same as noted above for $A^n_{k,l}$.

The main step in proving Theorem \pr{exact} is the following:

\begin{prop}\label{ankl}
\begin{enumerate}
\item For any $n\in \Z, k \in \zo$, $B^n_{k,k}$ is a basis for $\X^n_{k,k}$.
\item For any $n\in \Z, k < l \in \zo$, $B^n_{k,l}$ is an independent set spanning a subspace
of codimension 1 in $\X^n_{k,l}$.
\end{enumerate}
\end{prop}

\begin{pf}
We first look at the following model space.
Let $\E$ be the vector space with basis all symbols $E_i$, $i \in \Z$.
Let $D$ be the following set of element in $\E$:
$$D = \{ E_{-i} +E_i \}_{i \geq 0} \cup \{ E_{-i} +E_{1+i} \}_{i \geq 0}$$
We claim that $D$ is a basis for $\E$. Indeed,
order the symbols $E_i$ as follows: 
$$E_0, \ \ E_1, \ \ E_{-1}, \ \ E_2, \ \ E_{-2}, \ \ E_3, \ \ E_{-3}, \ \ \dots.$$ 
Order the elements of
$D$ by alternatingly taking an element from the first and second set, that is
$$2E_0, \ \ E_0 + E_1, \ \ E_{-1} + E_1, \ \ E_{-1} + E_2, \ \ E_{-2}+E_2 , \ \ E_{-2}+E_3, \ \ E_{-3}+E_3, \ \ \dots.$$
By induction, the span of the first $m$ elements in the first list coincides with that of the second list.  
This establishes our claim, and we now prove (1) and (2) of the proposition.

\emph{Proof of (1):} Assume first that $n$ is even, and let $n=2m$.
Define an isomorphism $\varphi : \E \to \X^n_{k,k}$ by $\varphi(E_i) = X^{m-i,m+i}_{k,k}$ for all $i \in \Z$. 
Then for all $i \geq 0$
\begin{itemize}
\item $\varphi(E_{-i} + E_i) = X^{m+i,m-i}_{k,k} + X^{m-i, m+i}_{k,k} = F^\1(J^+_{{ m+i \ch k }{ m-i \ch k }})$
\item $\varphi(E_{-i} + E_{1+i}) = X^{m+i,m-i}_{k,k} + X^{m-1-i, m+1+i}_{k,k}
= X^{m+i,(m-1-i)+1}_{k,k} + X^{m-1-i, (m+i)+1}_{k,k} = F^\1( J^A_{{m+i \ch k  }{ m-1-i \ch k }} )$ 
\end{itemize}
That is, $\varphi$ maps $D$ bijectively onto
$B^n_{k,k}$, and so $B^n_{k,k}$ is a basis for $\X^n_{k,k}$.

For $n$ odd, let $n=2m+1$. This time we define the isomorphism  
$\varphi : \E \to \X^n_{k,k}$ by $\varphi(E_i) = X^{m+i,m+1-i}_{k,k}$ for all $i \in \Z$.
We have for all $i \geq 0$,
\begin{itemize}
\item $\varphi(E_{-i} + E_i) = X^{m-i,(m+i)+1}_{k,k} + X^{m+i, (m-i)+1}_{k,k} = F^\1(J^A_{ { m-i \ch k }{ m+i \ch k } })$
\item $\varphi(E_{-i} + E_{1+i}) = X^{m-i,m+1+i}_{k,k} + X^{m+1+i, m-i}_{k,k} = F^\1(J^+_{{ m-i \ch k }{ m+1+i \ch k }})$. 
\end{itemize}
That is, again, $\varphi$ maps $D$ onto $B^n_{k,k}$, and so $B^n_{k,k}$ is a basis for $\X^n_{k,k}$.

\emph{Proof of (2):} Let $\E'$ be another copy of $\E$ defined with symbols $\{E'_i\}$ and construct the corresponding 
basis $D'$, then $D \cup D'$ is a basis for $\E \oplus \E'$. 
Assume first that $n$ is even and $n=2m$.
Let $\varphi:\E \oplus \E' \to \X^n_{k,l}$ be the isomorphism defined as follows:
\begin{itemize}
\item For $i \leq 0$: $\varphi(E_i) = X^{m-i,m+i}_{k,l}$ and for $i \geq 1$: $\varphi(E_i) = X^{m-i,m+i}_{l,k}$
\item For $i \leq 0$: $\varphi(E'_i) = X^{m-i,m+i}_{l,k}$ and for $i \geq 1$: $\varphi(E'_i) = X^{m-i,m+i}_{k,l}$
\end{itemize}
Then the value of $\varphi$ on $D \cup D'$ is as follows:
\begin{itemize}
\item $\varphi(2E_0) = 2X^{m,m}_{k,l} \not\in B^n_{k,l}$ 
\item $\varphi(2E'_0) = 2X^{m,m}_{l,k} \not\in B^n_{k,l}$
\end{itemize}

For $i \geq 1$:
\begin{itemize}
\item $\varphi(E_{-i} + E_i) = X^{m+i,m-i}_{k,l} + X^{m-i, m+i}_{l,k} = F^\1(J^+_{{ m+i \ch k },{ m-i \ch l }})$
\item $\varphi(E'_{-i} + E'_i) = X^{m+i,m-i}_{l,k} + X^{m-i, m+i}_{k,l} = F^\1(J^+_{{ m-i \ch k },{ m+i \ch l }})$
\end{itemize}
(Note that $F^\1( J^+_{ { m \ch k },{ m \ch l  } }  )$ is skipped here.)

For all $i \geq 0$:
\begin{itemize}
\item $\varphi(E_{-i} + E_{1+i}) = X^{m+i,(m-1-i)+1}_{k,l} + X^{m-1-i, (m+i)+1}_{l,k} 
         = F^\1(J^A_{{ m+i \ch k },{ m-1-i \ch l }})$
\item $\varphi(E'_{-i} + E'_{1+i}) = X^{m+i,(m-1-i)+1}_{l,k} + X^{m-1-i, (m+i)+1}_{k,l} 
         = F^\1(J^A_{{ m-1-i \ch k },{ m+i \ch l }})$
\end{itemize}

That is, the set $\varphi(D \cup D')$, which is a basis for $\X^n_{k,l}$, 
does not precisely coincide with $B^n_{k,l}$. The only difference is that the two elements 
$2X^{m,m}_{k,l}, 2X^{m,m}_{l,k}$ that appear in $\varphi(D \cup D')$, are replaced in 
$B^n_{k,l}$ by the one element $X^{m,m}_{k,l} + X^{m,m}_{l,k} = F^\1( J^+_{ { m \ch k },{ m \ch l  } }  )$. 
It follows that $B^n_{k,l}$ is independent and $spanB^n_{k,l}$ is 
of codimension 1 in $\X^n_{k,l}$.

For $n$ odd let $n=2m+1$.    
Let $\varphi:\E \oplus \E' \to \X^n_{k,l}$ be the isomorphism defined as follows:
\begin{itemize}
\item For $i \leq 0$, $\varphi(E_i) = X^{m+i,m+1-i}_{k,l}$ and for $i \geq 1$: $\varphi(E_i) = X^{m+i,m+1-i}_{l,k}$
\item For $i \leq 0$, $\varphi(E'_i) = X^{m+i,m+1-i}_{l,k}$ and for $i \geq 1$: $\varphi(E'_i) = X^{m+i,m+1-i}_{k,l}$
\end{itemize}
Then the value of $\varphi$ on $D \cup D'$ is as follows: 
\begin{itemize}
\item $\varphi(2E_0) = 2X^{m,m+1}_{k,l} \not\in B^n_{k,l}$ 
\item $\varphi(2E'_0) = 2X^{m,m+1}_{l,k} \not\in B^n_{k,l}$
\end{itemize}

For $i \geq 1$:
\begin{itemize}
\item $\varphi(E_{-i} + E_i) = X^{m-i,(m+i)+1}_{k,l} + X^{m+i, (m-i)+1}_{l,k} = F^\1(J^A_{{ m-i \ch k },{ m+i \ch l }})$
\item $\varphi(E'_{-i} + E'_i) = X^{m-i,(m+i)+1}_{l,k} + X^{m+i, (m-i)+1}_{k,l} = F^\1(J^A_{{ m+i \ch k },{ m-i \ch l }})$
\end{itemize}
(Note that $F^\1( J^A_{ { m \ch k },{ m \ch l  } }  )$ is skipped here.)

For all $i \geq 0$:
\begin{itemize}
\item $\varphi(E_{-i} + E_{1+i}) = X^{m-i,m+1+i}_{k,l} + X^{m+1+i, m-i}_{l,k} 
= F^\1(J^+_{{ m-i \ch k },{ m+1+i \ch l }})$
\item $\varphi(E'_{-i} + E'_{1+i}) = X^{m-i,m+1+i}_{l,k} + X^{m+1+i, m-i}_{k,l} 
= F^\1(J^+_{{ m+1+i \ch k },{ m-i \ch l }})$
\end{itemize}

Again, the set $\varphi(D \cup D')$, which is a basis for $\X^n_{k,l}$, 
does not coincide with $B^n_{k,l}$. The two elements 
$2X^{m,m+1}_{k,l}, 2X^{m,m+1}_{l,k}$ that appear in $\varphi(D \cup D')$, are replaced in $B^n_{k,l}$ 
by the one element $X^{m,m+1}_{k,l} + X^{m,m+1}_{l,k} = F^\1( J^A_{{m \ch k },{m \ch l}})$. Again 
it follows that $B^n_{k,l}$ is independent and $spanB^n_{k,l}$ is 
of codimension 1 in $\X^n_{k,l}$.
\end{pf}

Proposition \pr{ankl} states in particular that for any $n \in \Z, k \leq l \in \zo$,
$F^\1(A^n_{k,l})$ is an independent set.
This implies that $A^n_{k,l}$ itself is independent, and so a basis for $\J^n_{k,l}$,
and that $F^\1|_{\J^n_{k,l}} : \J^n_{k,l} \to \X^n_{k,l}$ is injective.
We have already noted that the symbols of type $J^+,J^A$ span $\J$
and so $\J$ is the sum of all $\J^n_{k,l}$. Since $\X$ is the direct sum of all $\X^n_{k,l}$,
and each $F^\1|_{\J^n_{k,l}}$ is injective, we obtain that $\J$ is the \emph{direct} sum of all $\J^n_{k,l}$.
(And since we have seen $A^n_{k,l}$ is a basis for $\J^n_{k,l}$, we also obtain 
that the set of all symbols of type $J^+,J^A$ is a basis for $\J$.)

For $n \in \Z$ and $k \leq l \in \zo$ define $\Y^n_{k,l} = span\{ Y^n_{k,l} \}\su \Y$.
Then for $k=l$ we have $dim \Y^n_{k,l}=0$ and for $k < l$ we have $dim \Y^n_{k,l} =1$.
$\Y$ is the direct sum of all $\Y^n_{k,l}$, and $\Psi(\X^n_{k,l}) \su \Y^n_{k,l}$.

Now, to prove Theorem \pr{exact} it is enough to show that for each $n \in \Z$ and $k \leq l \in \zo$,
the restricted sequence $0 \to \J^n_{k,l} @>F^\1>> \X^n_{k,l} @>\Psi>> \Y^n_{k,l} \to 0$ is exact.
Clearly $\Psi \circ F^\1 |_\J = 0$, and  $\Psi|_{\X^n_{k,l}} : \X^n_{k,l} \to \Y^n_{k,l}$ is surjective,
and we have also established that $F^\1|_{\J^n_{k,l}} : \J^n_{k,l} \to \X^n_{k,l}$ is injective.
Since $F^\1(\J^n_{k,l}) = span B^n_{k,l}$, and since by Proposition \pr{ankl}
$codim ( span B^n_{k,l} ) = dim \Y^n_{k,l}$, indeed 
$0 \to \J^n_{k,l} @>F^\1>> \X^n_{k,l} @>\Psi>> \Y^n_{k,l} \to 0$ is exact, and so
$0 \to \J @>F^\1|_{\J}>> \X @>\Psi>> \Y \to 0$ is exact.

\end{document}